\documentclass[preprint]{elsarticle}

\usepackage{hyperref,amsmath,amsfonts,amsthm}
\usepackage{dsfont}

\journal{Statistics and Probability Letters}

\newcommand{\cB}{{\mathcal B}}

\newcommand{\cS}{{\mathcal S}}

\newcommand{\GS}{{\mathds{S}}}

\newcommand{\GR}{{\mathds{R}}}

\newcommand{\bP}{{\mathbb{P}}}

\newcommand{\Idn}{\mbox{ $1\mskip-1.6\thinmuskip$I} }

\newcommand{\conver}{\mathop{\longrightarrow}}

\theoremstyle{plain}
\newtheorem{theorem}{\bf Theorem}

\theoremstyle{definition}

\newtheorem{remark}{\bf Remark}

\begin{document}

\begin{frontmatter}

\title{Quenched phantom distribution functions\\ for Markov chains}

\author[Torun]{Adam Jakubowski\corref{corresp}}
\cortext[corresp]{Corresponding author}
\ead{adjakubo@mat.umk.pl}

\author[Torun]{Patryk Truszczy\' nski}
\ead{502121@doktorant.umk.pl}

\address[Torun]{Nicolaus Copernicus University, Faculty of Mathematics and Computer Science, ul.~Chopina 12/18, 87-100 Toru\'n, Poland}





\begin{abstract}
It is known that random walk Metropolis algorithms with heavy-tailed target densities can model atypical (slow) growth of maxima, which in general is exhibited by processes with the extremal index zero. The asymptotics of maxima of such sequences can be analyzed in terms of continuous phantom distribution functions. We show that in a large class of positive Harris recurrent Markov chains (containing the above Metropolis chains) a phantom distribution function can be recovered by starting ``at the point" rather than from the stationary distribution.   
\end{abstract}


\end{frontmatter}


\section{Introduction and the statement of results} 
A stationary sequence $\{X_n\}$ of random variables with marginal distribution function $F(x) = \bP\big( X_j \leq x\big)$ and partial maxima $M_n = \max_{0 \leq j \leq n - 1} X_j$ admits a phantom distribution function $G$ if 
\begin{equation}\label{eq1}
\sup_{x\in\GR} \big| \bP\big(M_n \leq x\big) - G^n(x)\big| \conver 0, \quad \text{as $n\to \infty$}.
\end{equation}
Clearly, $G$ is not uniquely determined, since for (\ref{eq1}) to hold 
only the behavior of $G(x)$ near the right end $G_* = \sup \{x\,;\, G(x) < 1\}$
 is of importance. It was observed in \cite[Proposition 1]{DJL15} that under
 natural conditions (regularity), if $H$ is any other phantom distribution function
 for $\{X_n\}$, it must be tail equivalent to $G$, i.e.
\[ \lim_{x \to G_*-} \frac{1 - H(x)}{1 - G(x)} = 1.\]
   
Phantom distribution functions were introduced by O'Brien in \cite{OBr87}, 
as an extension of the notion of Leadbetter's extremal index. The latter
 corresponds to $G(x)$ of the form $F^{\theta}(x)$, for some $\theta \in (0,1]$
 (see \cite{Lead83}).  Contrary to extremal indices, the existence of phantom
 distribution functions is a quite common phenomenon for weakly dependent
 stationary sequences. For example, phantom distribution functions can be
 almost explicitly constructed for Markov chains with regenerative structure
 \cite[Theorem 3.1]{Root88}. Moreover, it was proved in the recent paper
 \cite{DJL15} that any $\alpha$-mixing stationary sequence with {\em
 continuous} marginals $F$ admits a continuous phantom distribution function.

It follows that the asymptotic behavior of maxima of stationary sequences 
with no extremal index or with the extremal index zero can still be analyzed
 using phantom distribution functions. Recall that $\{X_n\}$ has the extremal
 index $\theta = 0$ if 
\[ \bP\big( M_n \leq u_n(\tau)\big) \conver 1,\]
for every sequence $u_n(\tau)$ such that $n\big( 1 - F(u_n(\tau))\big) \to \tau$
 (see  \cite{Lead83}).  This means that maxima $M_n$ increase slower than
 maxima of i.i.d. random variables with marginal's tails comparable to 
$1-F(x)$, and so  asymptotic properties of $M_n$ cannot be expressed in
 terms of $F$.
Such a situation appears, for example, when Lindley's process has
 subexponential innovations (see \cite{Asmu98}) or when the continuous 
target distribution of the random walk Metropolis algorithm has heavy tails
 (see \cite{RRSS06}).
  
We will focus on the latter example. Let us recall basic definitions.
Let $\{Z_j\}$ be an i.i.d. sequence with the marginal distribution function $H$ given by the {\em proposal} density $h$, which is symmetric about $0$, and let $\{U_j\}$ be an i.i.d. sequence distributed uniformly
on $[0,1]$, independent of $\{Z_j\}$. Choose and fix the {\em target} probability density $f(x)$ (heavy-tailed in our case).
  Then the random walk Metropolis algorithm is the  Markov chain given by the recursive equation
\begin{equation}\label{eq:MH}
X_{j+1} = X_j + Z_{j+1} \Idn\big\{ U_{j+1} \leq \psi\big(X_j, X_j + Z_{j+1}\big)\big\},
\end{equation}
where  $\psi(x,y)$ is defined as
\begin{equation}\label{eq:MH2}
\psi(x,y) = \begin{cases}
\min\big\{f(y)/f(x), 1\big\} &\text{ if } f(x) > 0,\\
1 & \text{ if } f(x) = 0.
\end{cases}
\end{equation}

Random walk Metropolis algorithms have
 been designed for the purposes of simulation and therefore can be efficiently
 used in modeling of atypically slow growth of partial maxima of stationary
 sequences.  As observed in \cite[Remark 7]{DJL15} such a flexible family of models allows classifying stationary processes by relation of having a relative extremal index in the sense of \cite{Jak91}. To be more precise, given a ``model" stationary sequence $\{X_n\}$ we can consider a class of stationary sequences $\{X_n'\}$ satisfying for some $\theta \in (0,+\infty)$  
\begin{equation} \label{eq2}
\sup_{x\in\GR} \big| \bP\big(M_n \leq x\big) - \bP^{\theta}\big(M_n' \leq x\big)\big| \conver 0, \quad \text{as $n\to \infty$},
\end{equation}
where $M_n' = \max_{1 \leq j \leq n-1} X_j'$. Writing $\{ X_n\}\sim_{\theta} \{X_n'\}$ if (\ref{eq2}) holds, we see that $\{ X_n\}\sim_{\theta} \{X_n'\}$ implies $\{ X_n'\}\sim_{1/\theta} \{X_n\}$ and that $\bP\big(M_n \leq v_n\big) \to \alpha \in (0,1)$ implies $\bP\big(M_n' \leq v_n\big) \to \alpha^{1/\theta}$,
what clearly relates asymptotic quantiles of both sequences.

For weakly dependent sequences relation (\ref{eq2}) means that there exists a phantom distribution function $G$ for $\{X_n\}$ such that 
$G^{1/\theta}$ is a phantom distribution function for $\{X_n'\}$ (see the proof of \cite[Theorem 1.5]{Jak91}). It follows that the idea of comparison to models with atypical growth of maxima requires efficient methods of finding phantom distribution functions for random walk Metropolis algorithms with heavy-tailed proposals.  We refer to the discussion in \cite{JaRo07} for peculiarities 
related to the choice of the proposal density of the Monte Carlo Markov Chain algorithms.

Here we want to address another issue: whether we can start our
random walk Metropolis algorithms ``at the point'' rather than at the target density. The question is not trivial, since starting ``at the point'' we loose stationarity. 

We adopt formula (\ref{eq1}) as the definition of a phantom distribution function of a non-stationary sequence (see \cite{Jak93}). For the terminology related to general space Markov chains $\{Y_n\}$ we refer to the classic source \cite{MeTw09}. In particular, $\bP_s(\cdot)$ means the probability conditional on $\{ Y_0 = s\}$, while  $ \bP_{\pi}(\cdot)$ means that $Y_0$ is distributed according to $\pi$.

\begin{theorem}\label{main}
Let $\{Y_n\}$ be a positive Harris and aperiodic chain taking values in
$\big(\GS,\cS)$ and with a stationary distribution $\pi$. Let  $f: (\GS,\cS\big) \to \big(\GR^1,\cB^1\big)$ be a measurable function.

Let us define 
\[ X_n = f(Y_n),\quad n = 0,1,2,\ldots, \quad M_n = \max_{0 \leq j \leq n -1} X_j, \ n = 1,2, \ldots.\]

If $\{X_n\}$ admits a continuous phantom distribution function $G$ under some initial distribution $\lambda$, i.e. if we have 
\begin{equation}\label{eqlambda}
 \sup_{x\in\GR^1} \Big|\bP_{\lambda}\big( M_n \leq x\big) - G^n(x)\Big| \to 0, \text{ as $n \to \infty$},
\end{equation}
then $G$ is also a continuous phantom distribution function for the stationary (under $\pi$) sequence $\{X_n\}$.

Conversely, if $\{X_n\}$
admits a continuous phantom distribution function $G$ under $\pi$, then 
there exists a set $\GS_0 \in \cS$ satisfying $\pi\big(\GS_0\big) = 0$ and such that relation (\ref{eqlambda}) holds for every initial distribution $\lambda$ with the property that $\lambda(\GS_0)=0$ 
 \end{theorem}

\begin{theorem}\label{thalpha}
In assumptions of Theorem \ref{main}, if $\pi\circ f^{-1}$ (i.e. the marginal law of $\{X_n\}$ under the stationary distribution $\pi$) is continuous and unbounded above, then there exists a {\em continuous} distribution function $G$ such that for each $s\in\GS$
\[ \sup_{x\in\GR^1} \Big|\bP_{s}\big( M_n \leq x\big) - G^n(x)\Big| \to 0, \text{ as $n \to \infty$}.\]
\end{theorem}

The proofs are given in the next section. We will need a complement to \cite[Corollary 5]{Jak93}, which might be of independent interest.

Recall that a distribution function $G$ is regular (in the sense of O'Brien), if
\begin{equation}
\label{e2obr}
G(G_*-)=1\quad\text{and}\quad\lim\limits_{x\to G_*-}\dfrac{1-G(x-)}{1-G(x)}=1.
\end{equation}
Due to an observation made long time ago by \cite{OBr74} (Theorem 2), $G$ is regular if, and only if, for some $\gamma \in (0,1)$ there exists a sequence $\{v_n = v_n(\gamma)\}$ such that 
\[G^n(v_n) \to \gamma.\] 
Notice that if $G$ is regular then the sequence $\{v_n(\gamma)\}$ exists for every $\gamma \in (0,1)$ and that $\{v_n(\gamma)\}$ can always be chosen non-decreasing.

\begin{theorem}\label{key}
Let $X_0, X_1, X_2, \ldots, $ be an {\em arbitrary} sequence of random variables with partial maxima $M_n = \max_{0 \leq j \leq n-1} X_j$.  Then the following conditions are equivalent.
\begin{description}
\item{\bf (i)} $\{X_j\}$ admits a {\em continuous} phantom distribution function.
\item{\bf (ii)} $\{X_j\}$ admits a {\em regular} phantom distribution function. \item{\bf (iii)}There exists $\beta > 0$ and a non-decreasing sequence of levels $\{v_n = v_n(\beta)\}$ such that
\begin{equation}\label{eqkey}
\bP\big( M_{[nt]} \leq v_n(\beta)\big) \to \exp (- \beta t),\quad t \in D,
\end{equation}
where $D \subset  \GR^+$ is dense.
\item{\bf (iv)} For every $\beta > 0$ there exists a non-decreasing sequence of levels $\{v_n = v_n(\beta)\}$ such that (\ref{eqkey}) holds for every $t > 0$
(i.e. for $D = \GR^+$).
\end{description}
\end{theorem}

\begin{remark}
If (\ref{eqkey}) holds for some dense subset $D$ of $\GR^+$, then it is satisfied {\em uniformly } in $t\in \GR^+$.
\end{remark}

\begin{remark}
We have suggested one possible motivation for considering {\em quenched} (i.e. started ``at the point'') phantom distribution function. But such a notion is  interesting by itself, as in the context of the central limit theorem 
\cite{DeLi03} or the functional central limit theorem  \cite{BPP16}. It must be stressed that passing from a ``usual'' limit theorem to its quenched form need not be automatic (see \cite{VoWo10}).  

According to our knowledge this paper is the first that addresses ``quenched'' questions in the  extreme value limit theory.
\end{remark}

\section{Proofs}
\subsection{Proof of Theorem \ref{key}}

Statement {\bf (iv)} trivially implies {\bf (iii)}, so let us assume that {\bf (iii)} holds, i.e.  (\ref{eqkey}) is satisfied for some $\beta > 0$ and some non-decreasing sequence of levels $\{v_n(\beta)\}$. Applying 
Corollary 5 in \cite{Jak93} we obtain a phantom distribution function $G$ for  $\{X_j\}$ given by the formula
\begin{equation}\label{phdf}
G(x) = \begin{cases} 0, & \text{ if $x < v_1(\beta)$,}\\
\exp( - \beta)^{1/n}, & \text{ if $v_n(\beta) \leq x < v_{n+1}(\beta)$,}\\
1, & \text{ if $x \geq \sup_n v_n(\beta)$.}
\end{cases}
\end{equation}
$G$ of the above form is regular, hence {\bf (ii)} holds. By \cite[Remark 1, p. 704]{DJL15} there exists a {\em continuous} distribution function $H$ such that
\[ \sup_{x\in\GR^1} \big| G^n(x) - H^n(x)\big| \to 0, \quad \text{ as $n\to\infty$.}\]
Clearly, $H(x)$ is a continuous phantom distribution function for $\{X_j\}$. So {\bf (i)} also holds.
(The interested reader may find an explicit formula for $H$ on p.  704 in \cite{DJL15}).  

It remains to prove {\bf (i)} $\Rightarrow$ {\bf (iv)}.
So let us suppose that $\{X_j\}$ admits a continuous phantom distribution function, i.e. (\ref{eq1}) holds for some continuous $G$. Choose $\beta > 0$ and define 
\[ v_n(\beta) = \inf\{x\,;\, G^n(x) = \exp( - \beta)\}.\]
(notice that such $v_n(\beta)$ exists by the continuity of $G$). 
Then for any $t > 0$
\[
G^{[nt]}(v_n(\beta)) = \exp( - \beta)^{[nt]/n} \to \exp( - t\beta), \text{ as $n\to\infty$,}
\] 
and
\[ \bP\big( M_{[nt]} \leq v_n(\beta)\big) = G^{[nt]}(v_n) + o(1) \to \exp( - t\beta), \text{ as $n\to\infty$.}\]

\subsection{Proof of Theorem \ref{main}}

Let $\{Y_j\}$ be an aperiodic positive Harris chain and let $\pi$ be its unique stationary initial distribution. Choose and fix an initial distribution $\lambda$. By \cite[Proposition 3.13, p. 205]{Asmu03} there exist a coupling of $\{X_n\}$ under $\pi$ and $\lambda$ with some a.s. finite coupling time $\tau$. Recall, that this means that on some 
probability space one can define two stochastic processes $\{Y_j'\}$ and 
$\{Y_j^{\prime\prime}\}$ such that $\{Y_j'\}$ has the same distribution as $\{Y_j\}$ under the initial distribution $\pi$, $\{Y_j^{\prime\prime}\}$ has the same distribution as $\{Y_j\}$ under the initial distribution $\lambda$ and 
\[ Y_t'(\omega) = Y_t^{\prime\prime}(\omega), \text{ whenever $t \geq \tau(\omega)$.}\]
Let $X_j' = f(Y_j')$ and $X_j^{\prime\prime} = f(Y_j^{\prime\prime})$ and let $\{M_n'\}$ and  $\{M_n^{\prime\prime}\}$ be the partial maximum processes for $\{X_j'\}$ and $\{X_j^{\prime\prime}\}$, respectively.

Suppose that (\ref{eqlambda}) holds for some continuous $G$. By Theorem \ref{key} there exists $\beta > 0$ and a non-decreasing sequence of levels $v_n \nearrow G_*$ such that in the present notation
\begin{equation}\label{eq:keykey} 
\bP_{\lambda}\big( M_{[nt]} \leq v_n\big) = \bP\big( M_{[nt]}^{\prime\prime} \leq v_n\big) \to \exp(-\beta t),\quad t > 0.
\end{equation}

 We claim that 
\begin{equation}\label{eq:good}
\bP_{\lambda} \big( X_n < G_*\big) = 
\bP\big( X_n^{\prime\prime} < G_*\big) = 1, \quad n = 0,1,2,\ldots.
\end{equation}
Indeed, if $\bP\big( X_{n_0}^{\prime\prime} \geq G_*\big) = \delta > 0$ for some $n_0$, then for $n$ such that $nt \geq n_0$ we have 
\[ \bP\big( M_{[nt]}^{\prime\prime} \leq v_n\big) \leq  \bP\big( X_{n_0}^{\prime\prime} < G_*\big) = 1 - \delta,\]
while by (\ref{eq:keykey})
\[ \bP\big( M_{[nt]}^{\prime\prime} \leq v_n\big) \to \exp(-\beta t) > 1 - \delta,\]
if $t < - \ln (1-\delta)/\beta$. 

So assume (\ref{eq:good}). By the convergence in total variation of marginals to the stationary distribution (see e.g. \cite[Theorem 13.3.3, p. 328]{MeTw09})
we have also for $m=0,1,2,\ldots$
\begin{equation}\label{eq:goodgood}
\bP\big( X_m' < G_*\big) = \bP_{\pi} \big( X_1 < G_*\big) = \lim_{n\to\infty}  \bP_{\lambda} \big( X_n < G_*\big) = 1.
\end{equation}
 Since the coupling time $\tau$ is a.s. finite, we have
\[ \bP\big( \max_{0 \leq j \leq \tau-1} X_j' < G_*\big) = \sum_{k=1}^{\infty} 
 \bP\big( \max_{0 \leq j \leq k-1} X_j' < G_*, \tau = k\big) =  \sum_{k=1}^{\infty} \bP\big( \tau = k\big) = 1,\]
hence for every $t > 0$
\begin{align*}
 1 \geq \bP\big( \max_{ 0 \leq j \leq ([nt] \wedge\tau-1)} &X_j' \leq v_n \big)  \geq \\
&\geq  \bP\big( \max_{ 0 \leq j \leq \tau-1} X_j' \leq v_n\big) \nearrow \bP\big( \max_{ 0 \leq j \leq \tau-1} X_j' < G_*\big) = 1.
\end{align*}
In a similar way we obtain that
$ \bP\big( \max_{0 \leq j \leq \tau-1} X_j^{\prime\prime} < G_*\big) = 1$
and for every $t > 0$
\[
\bP\big( \max_{ 0 \leq j \leq ([nt] \wedge\tau-1)} X_j^{\prime\prime}\leq v_n \big)  \to 1.\]
Therefore we can write
\begin{align}
\bP_{\pi}\big( M_{[nt]} \leq v_n\big) &= \bP\big( M_{[nt]}' \leq v_n\big) \nonumber\\
&=  \bP\big(  \max_{ 0 \leq j \leq ([nt] \wedge\tau-1)} X_j'  \leq v_n, \max_{ [nt] \wedge\tau \leq j \leq [nt]-1} X_j'  \leq v_n\big)\nonumber \\
& = \bP\big(\max_{ [nt] \wedge\tau \leq j \leq [nt]-1} X_j'  \leq v_n\big) + o(1) \label{eq:align}\\
& = \bP\big(\max_{ [nt] \wedge\tau \leq j \leq [nt]-1} X_j^{\prime\prime} \leq v_n\big) + o(1) \nonumber \\
&=  \bP\big(  \max_{ 0 \leq j \leq ([nt] \wedge\tau-1)} X_j^{\prime\prime}  \leq v_n, \max_{ [nt] \wedge\tau \leq j \leq [nt]-1} X_j^{\prime\prime}  \leq v_n\big) + o(1) \nonumber\\
&=\bP\big( M_{[nt]}^{\prime\prime} \leq v_n\big) + o(1)  \to \exp(-\beta t),\quad t > 0. \nonumber
\end{align}
Applying Theorem \ref{key} we obtain that $G$ is a phantom distribution function for $\{X_n\}$ under $\pi$.

To prove the other part of Theorem 1 let us suppose $\{X_n\}$ admits a continuous phantom distribution function $G$ under $\pi$. Then by Theorem \ref{key}
\[\bP_{\pi}\big( M_{[nt]} \leq v_n\big) = \bP\big( M_{[nt]}' \leq v_n\big) \to \exp(-\beta t),\quad t > 0,\]
for some $\beta > 0$ and a non-decreasing sequence of levels
$v_n \nearrow G_*$. Similarly as for (\ref{eq:good}) we deduce that 
\[ \bP_{\pi} \big( X_1 < G_*\big) = 
\bP\big( X_n' < G_*\big) = 1, \quad n = 0,1,2,\ldots.\]

Let us define ``the bad set''
\begin{equation}\label{eq:bad}
 B_0= \{ y\,;\, f(y) \geq G_*\}.
\end{equation}
Since $\pi (B_0) = 0$ and the chain $\{X_n\}$ is $\pi$-irreducible (see \cite[Theorem 10.4.9, p. 246]{MeTw09}) we obtain that also the set of states that communicate with  $B_0$ has $\pi$-measure zero (see \cite[Proposition 4.2.2, p. 83]{MeTw09}):
\[ \pi\big( \{ y\,;\, P^n(y, B_0) > 0 \ \text{ for some $n \geq 1$}\}\big) = 0.\]
So let us set 
\[ \GS_0 = B_0 \cup \{ y\,;\, P^n(y, B_0) > 0 \ \text{ for some $n \geq 1$}\},\]
and assume that an initial distribution $\lambda$ has the property that
\[ \lambda\big(\GS_0\big) = 0.\]
Then we have 
\[ \bP_{\lambda} \big( X_n < G_*\big) =  \bP\big( X_n^{\prime\prime} < G_*\big) = 1, \quad n = 0,1,2,\ldots.\]

This means that the two crucial relations (\ref{eq:good}) and (\ref{eq:goodgood}) are satisfied and so we may repeat step by step the reasoning in (\ref{eq:align}), but in reverse order.
The theorem follows. 

\subsection{Proof of Theorem \ref{thalpha}}

By virtue of \cite[Theorem 13.3.3]{MeTw09} one obtains $\alpha$-mixing of 
$\{X_n\}$ under stationary distribution $\pi$ (for more details see e.g. \cite{Brad07} or \cite{Douk94}). 
  
If $\pi\circ f^{-1}$ is continuous then $\{X_n\}$ admits a continuous phantom distribution function $G$ by \cite[Theorem 6]{DJL15}. Moreover, since $\pi\circ f^{-1}$ is unbounded above, $G_* = +\infty$ and the bad set $B_0$ given by (\ref{eq:bad}) is empty. Hence our Theorem \ref{main} states that $G$ is a phantom distribution function for $\{X_n\}$ under any initial distribution $\lambda$. In particular for every $s\in\GS$ 
\[ \sup_{x\in\GR^1} \Big|\bP_{s}\big( M_n \leq x\big) - G^n(x)\Big| \to 0, \text{ as $n \to \infty$}.\]



\begin{thebibliography}{}

\bibitem[Asmussen, 1998]{Asmu98}
Asmussen, S. (1998).
\newblock Subexponential asymptotics for stochastic processes: extremal
  behavior, stationary distributions and first passage probabilities.
\newblock {\em Ann. Appl. Probab.}, 8:354--374.

\bibitem[Asmussen, 2003]{Asmu03}
Asmussen, S. (2003).
\newblock {\em \sc Applied Probability and Queues}.
\newblock Springer, New York.

\bibitem[Barrera et~al., 2016]{BPP16}
Barrera, D., Peligrad, C., and Peligrad, M. (2016).
\newblock Macroscopic models for long-range dependent network traffic.
\newblock {\em Stochastic Process. Appl.}, 126:1885--1900.

\bibitem[Bradley, 2007]{Brad07}
Bradley, R. (2007).
\newblock {\em \sc Introduction to strong mixing conditions}, volume~2.
\newblock Kendrick Press, Heber City UT.

\bibitem[Derriennic and Lin, 2003]{DeLi03}
Derriennic, Y. and Lin, M. (2003).
\newblock The central limit theorem for markov chains started at a point.
\newblock {\em Probab. Theory Relat. Fields}, 125:73--76.

\bibitem[Doukhan, 1994]{Douk94}
Doukhan, P. (1994).
\newblock {\em \sc Mixing: Properties and Examples}, volume~85 of {\em Lecture
  Notes in Statist.}
\newblock Springer, New York.

\bibitem[Doukhan et~al., 2015]{DJL15}
Doukhan, P., Jakubowski, A., and Lang, G. (2015).
\newblock Phantom distribution functions for~some stationary sequences.
\newblock {\em Extremes}, 18:697--725.

\bibitem[Jakubowski, 1991]{Jak91}
Jakubowski, A. (1991).
\newblock Relative extremal index of~two stationary processes.
\newblock {\em Stochastic Process. Appl.}, 37:281--297.

\bibitem[Jakubowski, 1993]{Jak93}
Jakubowski, A. (1993).
\newblock An asymptotic independent representation in~limit theorems for~maxima
  of~nonstationary random sequences.
\newblock {\em Ann. Probab.}, 21:819--830.

\bibitem[Jarner and Roberts, 2007]{JaRo07}
Jarner, S. and Roberts, G. (2007).
\newblock Convergence of~heavy-tailed Monte Carlo Markov chain algorithms.
\newblock {\em Scand. J. Stat.}, 34:781--815.

\bibitem[Leadbetter, 1983]{Lead83}
Leadbetter, M.~R. (1983).
\newblock Extremes and local dependence in~stationary sequences.
\newblock {\em Z. Wahrscheinlichkeitstheorie verw. Gebiete}, 65:291--306.

\bibitem[Meyn and Tweedie, 2009]{MeTw09}
Meyn, S. and Tweedie, R. (2009).
\newblock {\em \sc Markov Chains and Stochastic Stability. 2nd Ed.}
\newblock Cambridge University Press, Cambridge.

\bibitem[O'Brien, 1974]{OBr74}
O'Brien, G.~L. (1974).
\newblock Limit theorems for the maximum term of a stationary process.
\newblock {\em Ann. Probab.}, 2:54--545.

\bibitem[O'Brien, 1987]{OBr87}
O'Brien, G.~L. (1987).
\newblock Extreme values for~stationary and Markov sequences.
\newblock {\em Ann. Probab.}, 15:281--292.

\bibitem[Roberts et~al., 2006]{RRSS06}
Roberts, G.~O., Rosenthal, J., Segers, J., and Sousa, B. (2006).
\newblock Extremal indices, geometric ergodicity of~Markov chains and MCMC.
\newblock {\em Extremes}, 9:213--229.

\bibitem[Rootz\'en, 1988]{Root88}
Rootz\'en, H. (1988).
\newblock Maxima and exceedances of stationary Markov chains.
\newblock {\em Adv. Appl. Probab.}, 20:371–390.

\bibitem[Voln\'y and Woodroofe, 2010]{VoWo10}
Voln\'y, D. and Woodroofe, M. (2010).
\newblock An example of non-quenched convergence in the conditional central
  limit theorem for partial sums of linear processes.
\newblock In Berkes, I., Bradley, R., Dehling, H., Peligrad, M., and Tichy, R.,
  editors, {\em \sc Dependence in Probability, Analysis and Number Theory},
  pages 317--322. Kendrick Press, Heber City UT.

\end{thebibliography}

\end{document}